\documentclass[12pt]{amsart}

\usepackage{amsmath,amssymb,latexsym}
\usepackage{amsfonts,amstext,amsthm,amscd}
\usepackage{graphicx}
\usepackage{amsthm}
\theoremstyle{plain}
\usepackage{color}
\usepackage[all]{xy}

\newtheorem*{tmaA}{Theorem A}
\newtheorem{theorem}{Theorem}[section]

\newtheorem{deff}[theorem]{Definition}
\newtheorem{ex}[theorem]{Example}
\newtheorem{lemma}[theorem]{Lemma}

\newtheorem{prop}[theorem]{Proposition}

\numberwithin{equation}{section}

\def\eproof{\noindent{\hfill $\blacksquare$}\bigskip}

\newcommand{\CC}{\mathbb{C}}
\newcommand{\ER}{\widehat{\mathbb{C}}}
\newcommand{\ZZ}{\mathbb{Z}}

\newcommand{\KK}{\mathcal{K}}

\newcommand{\Fat}{\mathcal{F}}
\newcommand{\Jul}{\mathcal{J}}

\newcommand{\ninf}{\rightarrow\infty}

\begin{document}
\title{Commuting functions in class $\KK$ and families of wandering Baker domains.}
\author{Adri\'an Esparza Amador}
\begin{abstract}
Given commuting functions $f,g$, with at most a countable compact set of essential singularities, recent results for entire functions are extended to prove that Julia sets match, $\Jul(f)=\Jul(g)$, in a particular case in the class $\KK$. With this result on hand, a general form of functions in this class, is constructed with \textit{wandering Baker domains}. 
\end{abstract}
\maketitle

\section{Introduction}
Given a rational or entire transcendental function $f$, the family of iterates $\{f^n\}_{n\geq0}$ generates a dichotomy over the Riemann sphere: the called \textit{Fatou and Julia sets}. The Fatou set of function $f$, denoted by $\Fat(f)$, is defined as the set of points where the family of iterates form a normal family in a neighborhood of such points. The Julia set is the complement over the Riemann sphere of the Fatou set, general theory on rational iteration can be found in the books \cite{bea}, \cite{cg} and \cite{mil1}, for reference on entire functions see \cite{sch}.\\
By definition, the Fatou set is open and the Julia set is closed. The Fatou set is decomposed in connected components, which can be classified in general terms as pre-periodic, periodic or wandering. It is well known that rational functions has no wandering domains \cite{sul1}. For entire functions several examples has been exhibited, some of the first examples are described in \cite{bak2}. \\
Two functions $f,g$ are called \textit{commuting} or \textit{permutable} if they satisfy $f(g)=g(f)$. For commuting functions, an immediate question arise: what is the relation of the Fatou and Julia sets of such functions? In the case of rational functions, Fatou in \cite{fa1} and Julia in \cite{ju1} proved that $\Jul(f)=\Jul(g)$. The proof is based on the compactness of the Riemann sphere, where rational functions are well defined and are continuous. Then, such proof can not be extend for entire transcendental functions, since they are defined in the complex plane, which is not compact. In fact, the possible presence of wandering and periodic Baker domains (see next section for formal definition) cause the problem, since in both Fatou components, the iterates tend to the essential singularity, the point at infinity. \\
In \cite{bak2}, Baker solve the problem for particular commuting entire transcendental functions. The theorem statement is as follows.
\begin{theorem}[\cite{bak2}, Lemma 4.5]\label{comm1}
Suppose that $f$ and $g$ are entire transcendental functions, $f$ commutes with $g$, and $f=g+c$, where $c$ is some constante. Then $\Jul(f)=\Jul(g)$. 
\end{theorem}
In the present work, the ideas in \cite{bak2} are extended for commuting meromorphic functions with a countable compact set of essential singularities, the so called \textit{class} $\KK$ (see section below for formal definition). 
\begin{tmaA}\label{tmaA}
Let $f,g\in\KK$ and suppose that $f$ commutes with $g$ and that $f=g+c$, for some constante $c$. Then $\Jul(f)=\Jul(g)$ and $\Fat(f)=\Fat(g)$. 
\end{tmaA}
This theorem allows us to construct \emph{wandering Baker domains}. 
\section{Class $\KK$}
A particular characteristic of rational and entire transcendental functions, is that they are closed under composition, that is, if $f$ and $g$ are rational (or entire) functions, then $f(g)$ and $g(f)$ are again rational (or entire) functions. \\
In a meromorphic transcendental function, not only the point at the infinity is an essential singularity, but such singularity has preimagenes, poles of the functions. In the case there is only one pole and also is an omitted value, this class of function is also closed under composition, sometimes called the Radstr\"om class (see \cite{rad} for reference). When the pole is not an omitted value, or the functions has more than one pole, the iterates are no longer closed under composition, in fact, this general class is not closed under composition. \\
Consider for example the function
$$f(z)=\frac{1}{z-z_0}+\exp(z),\ z_0\in\CC,$$
which is meromorphic with $z_0$ as its pole. If we consider the second iterate
$$f^2(z)=f(f(z))=\frac{1}{\frac{1}{z-z_0}+\exp(z)-z_0}+\exp\left(\frac{1}{z-z_0}+\exp(z)\right),$$
it can be notice now that $z_0$ is also an essential singularity of $f^2$ as well as the point infinity $\infty$. In the same way, $f^3$ will have infinitely many essential singularities (each solution of the equation $\frac{1}{z-z_0}+\exp(z)-z_0=0$). Then, each iterate no longer belongs to the meromorphic class of functions.\\
To solve this problem, it is possible to define a more general class of meromorphic functions, which includes the meromorphic transcendental functions, that are closed under composition. We then define, the class of general transcendental meromorphic functions as follows.
\begin{deff}\label{fordef}
We say that a function $f$ belongs to class $\KK$ of \textbf{general meromorphic functions} if there exists a countable compact set $A(f)\subset\ER$ such that $f$ is an analytic function in $\ER\backslash A(f)$ but in no other superset. We say that $f$ is \textbf{general transcendental meromorphic function} if $A(f)\neq\emptyset$. 
\end{deff} 
The class $\KK$ is closed under composition. 
\begin{prop}[\cite{bo}, Proposition 1.3]
If $f\in\KK$, then  $f^2\in\KK$ with $A(f^2)=A_2(f)=A(f)\cup f^{-1}(A(f))$. 
\end{prop}
Fundamentals, basic concepts and results on iteration of functions in class $\KK$, can be found in the dissertations of A. Bolsch \cite{bo} and M. Herring \cite{her}, and the recent papers \cite{bdh} and \cite{dms}. Among many results, the dynamical classification of the Julia set for rational functions, holds for class $\KK$. This result will be useful in the proof of Theorem A. \\
For a function $f\in\KK$, a point $z_0\in\ER$ is a \textit{periodic point} of $f$ of period $p>0$ if $f^p(z_0)=z_0$ and $f^j(z_0)\neq z_0$ for $0<j<p$. We classified $z_0$ as
\begin{itemize}
\item \emph{attracting} if $|(f^p)'(z_0)|<1$,
\item \emph{repelling} if $|(f^p)'(z_0)|>1$, 
\item \emph{indifferent} if $|(f^p)'(z_0)|=1$.
\end{itemize}
For an indifferent periodic point $z_0$, if $(f^p)'(z_0)=1$, $z_0$ is also called \emph{parabolic}. When $p=1$, a periodic point is called a \emph{fixed point}. 
\begin{lemma}[\cite{bo}, Theorem 1.12]\label{repel}
If $f\in\KK$, the $\Jul(f)$ is the closure of its repelling periodic points.
\end{lemma}
\subsection{Families of Baker domains}
We focus our attention in two kinds of Fatou components, wandering domains and Baker domains defined as follows.
\begin{deff}
Let $U\subset\Fat(f)$ be a maximal domain of normality so that $f:U\rightarrow U$ is analytic. If there exist $z_0\in\partial U$ such that $f^n(z)\to z_0$ as $n\ninf$ for $z\in U$, and $f(z_0)$ is not defined, then $U$ is called an invariant \emph{Baker domain} of $f$ and $z_0$ is called its \emph{Baker point}.
\end{deff}
Note that, by definition, Baker domains does not exist for rational functions, and $z_0\in A(f)$ in definition \ref{fordef}. In \cite{dms} authors called the point $z_0$ an \emph{absorbing point}. The first example of Baker domain was given by Fatou in \cite{fa2}, proving that for the entire function $f(z)=z+1+e^{-z}$ there exists a half plane where $f^n\rightarrow\infty$. Basic properties of these Fatou components can be found in the survey \cite{r}. \\
In \cite{tes} author uses the ideas in \cite{rs} to provide sufficient conditions over a large class of functions in class $\KK$ that exhibit an infinite collection of invariant Baker domains. The following result is then obtained in \cite{tes}. 
\begin{theorem}\label{fam}
Let $\Omega\in\{\CC,\ER\}$ and $g:\Omega\rightarrow\ER$ be a transcendental meromorphic or rational function. Consider a function in class $\KK$ given by 
$$f(z)=z+\emph{exp}(g(z)),$$ 
with $A(f)=\overline{g^{-1}(\infty)}$. If $z_0\in g^{-1}(\infty)$ is a pole of $g$ of order $p\geq1$, then $f(z)$ has $p$-families of Baker domains with $z_0$ as its Baker point. Each family lies in a sector of opening $2\pi/p$ in a small neighborhood of $z_0$.  \\
\end{theorem}
Theorem \ref{fam}, will help us to construct functions in class $\KK$ with \emph{wandering Baker domains}. 
\subsection{Wandering domains}
As was mentioned in the introduction, Sullivan in \cite{sul1} proved that rational functions has no wandering domains. The existence of wandering domains for entire (or meromorphic) transcendental functions, was considered at the end of the last century. Authors as Baker (\cite{bak1}, \cite{bak2}), and Herman (\cite{hrm}) exhibit some constructions for wandering domains for entire functions. The ideas in the construction of Herman are described below, see also \cite{bak2}.\\
Consider an entire function 
$$f(z)=z+\varphi(z),$$
where $\varphi$ is entire and periodic of period $2\pi i$. Let $\pi:\CC\rightarrow\CC^*$ be the covering map given by $\pi(z)=e^{\pm z}$. Construct then the following diagram.

\begin{displaymath}
\xymatrix{
\CC \ar[r]^f \ar[d]_{\pi} & \CC\ar[d]^{\pi}\\
\CC^* \ar[r]^{g} & \CC^*} 
\end{displaymath}
And let $f$ be such that $g$ is an entire function also. Then, thanks to Bergweiler's results on lifting of Julia sets \cite{berg}, $g$ may transfer a description of the dynamic of $h^n$ to that of $f^n$ and vice versa. With this construction on hand, we consider the following example.\\

\textbf{Example:} (Example 5.1 in \cite{bak2}) Take the function 
$$f(z)=z-1+e^{-z},$$ 
and the covering map $\pi(z)=e^{-z}$. This way, we obtain the function
$$g(t)=ete^{-t},$$
which has a \emph{super-attracting fixed point} at $t=1$.  This point lift to infinitely many fixed super-attracting fixed points of $f$, $z_n=2n\pi i$, with mutually disjoints domains of attraction $D_n$. Since
$$f^n(z+2\pi i)=f^n(z)+2\pi i,$$
we have $\Fat(f)$ is translated into itself by $z\mapsto z+2\pi i$. \\
Set
$$g(z)=2\pi i+f(z)=z-1+e^{-z}+2\pi i,$$ 
so that $f(g)=g(f)$, then by Lemma \ref{comm1}, $\Jul(f)=\Jul(g)$ and $\Fat(f)=\Fat(g)$. Thus $D_n$ is a component of $\Fat(g)$ and $g(D_n)=2\pi i+f(D_n)=D_{n+1}$, and then each component $D_n$ is a wandering domains for $g$. 
\begin{figure}
\includegraphics [width=10cm]{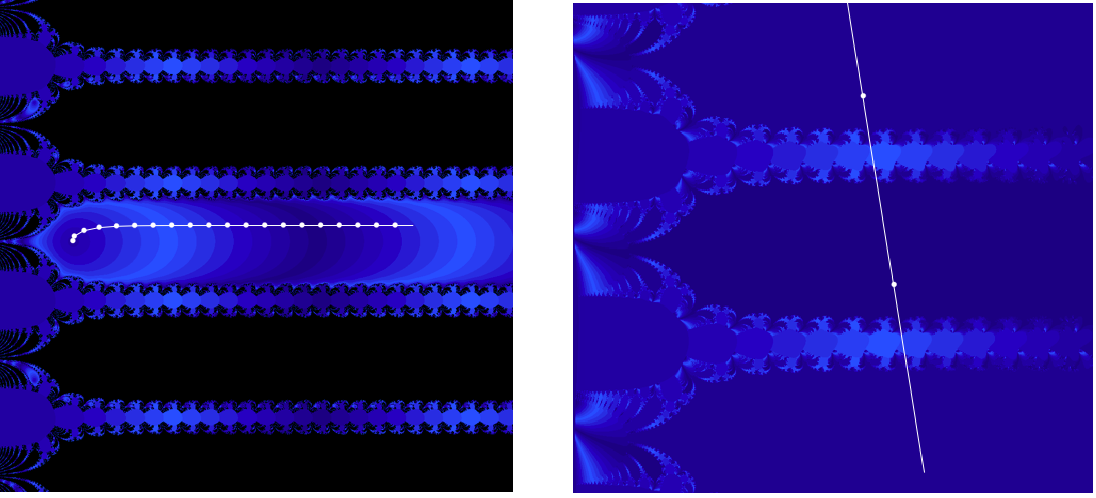}
\caption{On the left, we have the super-attracting domains in each band of width $2\pi$ for $f(z)=z-1+e^{-z}$. On the right, those super-attracting domains become wandering domains for $f(z)=z-1+e^{-z}+2\pi i$.}
\end{figure}
\subsection{Wandering Baker domains}
As in the previous example, consider now the entire function
$$f(z)=z+e^{-z}.$$
In \cite{bd}, the authors proved that $f$ contains a family $\{B_n\}_{n\in\ZZ}$ of Baker domains associated to the point at infinity (also, $f$ trivially holds conditions in Theorem 1 of \cite{rs}). Following the same construction, we obtain the entire function
$$h(t)=te^{-t},$$
with a \emph{parabolic fixed point} at $t=0$. This point is lifted to the point at infinity as the Baker point for each of the Baker domains $B_n$ of the function $f$. Also $\Fat(f)$ is translated into itself by $z\mapsto z+2\pi i$.  \\
Take $g(z)=2\pi i+f(z)$, so that $f(g)=g(f)$, then again by Lemma \ref{comm1}, $\Jul(f)=\Jul(g)$ and $\Fat(f)=\Fat(g)$. Thus $B_n$ is component of $\Fat(g)$ and $g(B_n)=2\pi i+f(B_n)=B_{n+1}$, and then each Baker domain $B_n$ of $f$ is a wandering domain for $g$. The fact that each $B_n$ is a Baker domain for $f$, makes that $B_n$'s are called \emph{wandering Baker domains} of $g$. 
\begin{figure}
\includegraphics [width=10cm]{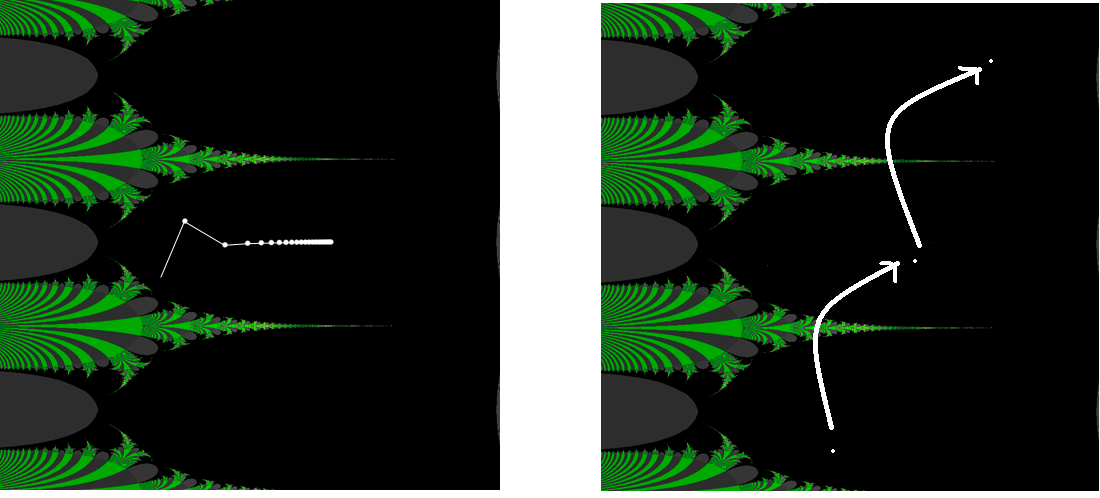}
\caption{On the left, there are Baker domains separated by the lines $y=(2k+1)\pi i$ for $f(z)=z+e^{-z}$. On the right, those Baker domains became in wandering domains vertically translated for $f(z)=z+e^{-z}+2\pi i$.}
\end{figure}
\section{Proof of Theorem A}
The proof is almost verbatim to the proof of Lemma 4.5 in \cite{bak2}. With the following variations. First, due the presence of several essential singularities different from $\infty$, compact sets should be considered at a positive distance away from $A(f)$, instead of just compact set of the complex plane. Second, Baker domains are allowed to be associated to finite points, then the proof has to be modified to finite Baker points. \\
\subsection{Proof} 
First, note that $A(f)=A(g)$ since $f=g+c$ with $c$-constant. Since Fatou and Julia sets are complementary in the Riemann sphere, and there is a symmetry of $f$ and $g$ about the hypothesis of the theorem, it is enough to prove:
\begin{itemize}
\item[a)] $g(\Jul(f))\subset\Jul(f)$ and
\item[b)] $g(\Fat(f))\subset\Fat(f)$.
\end{itemize}
To prove $a)$, from Lemma \ref{repel}, we have that $z_0\in\Jul(f)$ if and only if $z_0$ is a limit point of repelling periodic points of $f$, which are in $\Jul(f)$. Then is enough to prove that $g$ maps repelling periodic points of $f$ in repelling periodic points of $f$. By hypothesis $f(g)=g(f)$ in $\ER\backslash A(f)=\ER\backslash A(g)$, then 
$$f^p(g(z))=g(f^p(z)),$$
where this is well defined. So, if $z_0$ is a repelling periodic point of $f$, computing the derivative in the above expression, we have 
$$(f^p)'(g(z_0))\cdot g'(z_0)=g'(z_0)\cdot (f^p)'(g(z_0)),$$
which yields
$$(f^p)'(g(z_0))=(f^p)'(z_0).$$
In the case $g'(z)=0$, that is, $z_0$ is a critical point of $g$, we know that such set of points is an isolated set (discrete), since $\Jul(f)$ is perfect, the results still holds by the compactness of $\Jul(f)$. $a)$ is then proved. \\
To prove $b)$, consider the following two cases:\\
\begin{itemize}
\item[CASE 1] $f^n(z)\nrightarrow e$ for every $e\in A(f)$, that is, $z$ is in a periodic component of $\Fat(f)$ which is not a Baker domain.  
\item[CASE 2] $f^n(z)\rightarrow e$ for some $e\in A(f)$, that is, $z$ is contained in a Baker domain or a wandering domain.
\end{itemize}

CASE 1: Let $z_0\in\Fat(f)$ and $U$ be a neighborhood of $z_0$ such that $\overline{U}\subset\Fat(f)$, since $f^{n_k}$ does not converge to a point in $e\in A(f)=A(g)$, all $f^{n_k}(U)$ are contained in a compact set $K\subset\ER$ away from $A(f)$, $\text{dist}(K,A(f))>0$ in the spherical distance. Hence $g$ is uniformly continuous in $K$. Choose $U$ small enough such that $g(f^{n_k}(U))=f^{n_k}(g(U))$ has small diameter for large $n_k$. Then, by Montel's theorem, $g(U)\subset\Fat(f)$, in particular $g(z_0)\in\Fat(f)$. The first case is covered. \\

CASE 2: Let $z_0\in\Fat(f)$ and $U$ a neighborhood of $z_0$ where $f^n\rightarrow e$ as $n\ninf$ for some $e\in A(f)$. Take $\epsilon>0$ such that $\epsilon>|c|+1$, where $c=f-g$. There exist $n_0$ such that $|f^n-e|<1$ holds in $U$ for each $n>n_0$, then $|f(z)-e|<1$ for every $z\in f^n(U)$, with $n>n_0$. \\
If $g(z_0)\notin\Fat(f)$, then $\{f^m\}$ takes every value with at most two exceptions, for $m$ large in $g(U)$. So, there exist $t=g(\xi)$, with $\xi\in U$, such that for some $m>n_0$
$$|f^m(t)-e|=|f^m(g(\xi))-e|=|g(f^m(\xi))-e|>\epsilon.$$
Hence, $\eta=f^m(\xi)\in f^m(U)$ with $|g(\eta)-e|>\epsilon$, and
$$|c|=|g(\eta)-f(\eta)|=|(g(\eta)-e)-(f(\eta)-e)|>\epsilon-1,$$
which contradicts the choice of $\epsilon$. Then $g(z_0)\in\Fat(f)$ and $b)$ and the theorem is proved.\eproof

\section{Examples of wandering Baker domains}	
In the present section, we exhibit some examples of \emph{wandering Baker domains} in the class $\KK$, using the construction in Section 2. \\
\subsection{Construction}
For a periodic meromorphic function $g$ of period $P$, given Theorem \ref{fam}, the function
$$f(z)=z+\exp(g(z)),$$
possesses families of Baker domains, each family with \emph{Baker point} as a pole of function $g$. \\
Consider function $h\in\KK$, given by $h=f+P$, where $P$ is the period of $g$, that is,
$$h(z)=f(z)+P=z+\exp(g(z))+P.$$
Then, from direct calculation, we obtain that 
\begin{displaymath}
\begin{array}{rcl}
f(h(z))&=&f(z+\exp(g(z))+P)\\
       &=&z+\exp(g(z))+P+\exp(g(z+\exp(g(z))+P))\\
       &=&z+\exp(g(z))+\exp(g(z+\exp(g(z))))+P\\
       &=&f(z)+\exp(g(f(z)))+P\\
       &=&h(f(z)).
\end{array}
\end{displaymath}
Then $h$ commutes with $f$, using Theorem A we have $\Jul(f)=\Jul(h)$ and $\Fat(f)=\Fat(h)$, moreover, $\Fat(f)$ maps into itself by the map $z\mapsto z+P$. Then, the families of Baker domains of $f$ are actually families of \emph{wandering Baker domains} for $h$. \\
\subsection{Examples}
\begin{ex}\label{memo}
Consider the periodic meromorphic function $g(z)=\sin (z)$ with period $P=2\pi$. In \cite{dms}, it is proved, using \emph{Tracks Theory}, that the function
$$f(z)=z+\exp\left(\frac{1}{\sin (z)}\right),$$
has families of Baker domains associated to each zero of Sine function (the essential singularities of $f$). And then, that
$$h(z)=z+\exp\left(\frac{1}{\sin (z)}\right)+2\pi,$$
has families of wandering Baker domains. 
\end{ex}
\begin{figure}
\includegraphics [width=10cm]{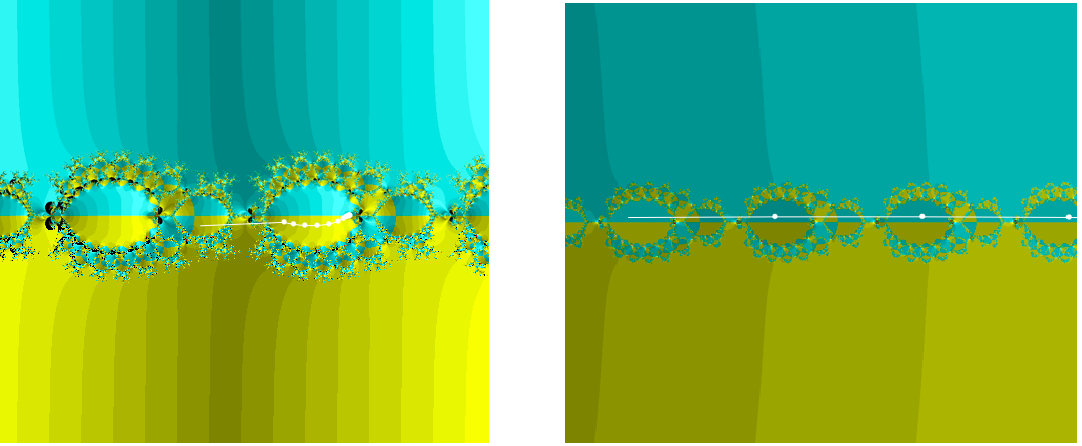}
\caption{On the left, finite Baker domains can be appreciated with Baker points over the real line for function $f(z)=z+\exp\left(\frac{1}{\sin (z)}\right)$. On the right, those Baker domains became wandering Baker domains for $h(z)=z+\exp\left(\frac{1}{\sin (z)}\right)+2\pi$.}
\end{figure}
\begin{ex}
Consider the periodic meromorphic function $g(z)=\frac{1}{1-\exp(z)}$ with period $2\pi i$, in the same way, by Theorem \ref{fam}, 
$$f(z)=z+\exp\left(\frac{1}{1-\exp(z)}\right)$$
has families of Baker domains associated to each solution of equation $1-\exp(z)=0$. Now, using Theorem A, 
$$g(z)=z+\exp\left(\frac{1}{1-\exp(z)}\right)+2\pi i$$
has wandering Baker domains.
\end{ex}

\textbf{Acknowledgements}. I want to thank G. Sienra for share with me his wonderful ideas in the subject and to encourage me to write this paper. In addition, I would like to thank the referee for his wise corrections.

\end{document}